\documentclass[12pt]{article}
\usepackage{amssymb}
\usepackage{amsmath}
\usepackage{latexsym}
\usepackage{epsfig}
\usepackage{color}

\newcommand {\eproof}{\space{\ \vbox{\hrule\hbox{\vrule height1.3ex\hskip0.8ex\vrule}\hrule}}\vskip 0.3cm \par}

\def\msn{\medskip\noindent}

\newcommand {\R}	{{\mathbf R}}

\newcommand {\Rn}       {\R^n}

\newcommand {\Rp}       {\R_+}
\newcommand {\Rpn}       {\R^n_+}

\newcommand {\Rmaxm}	{\R_{\max,\,\text{m}}}

\newcommand {\Rmaxmn}	{\R_{\max,\,\text{m}}^n} 

\newcommand {\ldotsn}	{{1,\ldots,n}}

\newcommand {\0}	{\mathbf{0}}

\newcommand {\1}	{\mathbf{1}}

\newcommand {\Hilb}	{\text{H}}

\newcommand {\cK}	{{\mathcal K}}

\newcommand {\cV}	{{\mathcal V}}

\newcommand {\odiv}	{/}

\newcommand {\leqoplus} {\leq_{\oplus}}

\newcommand {\leqoplusV}{\leq_{\oplus,\cV}} 

\newcommand {\supp}	{\operatorname{supp}}

\newcommand {\aff}	{\operatorname{aff}}

\newtheorem{theorem}          {Theorem}
\newtheorem{lemma}         [theorem]{Lemma}
\newtheorem{definition}    [theorem]{Definition}

\newtheorem{proposition}     [theorem]{Proposition}

\begin{document}

\title{Cyclic projectors and separation theorems
in idempotent convex geometry\thanks{%
Supported by the RFBR grant 05-01-00824 and the joint
RFBR/CNRS grant 05-01-02807}}
\author{St{\'{e}}phane Gaubert\thanks%
{INRIA, Rocquencourt, B.P. 105, 78105 Le Chesnay cedex, France. 
E-mail: Stephane.Gaubert@inria.fr} and {Serge{\u{\i}} Sergeev\thanks%
{Department of Physics, Sub-Department of Quantum Statistics and 
Field Theory, Moscow State University, Moscow,
119992 Leninskie Gory, Russia.
E-mail: sergiej@gmail.com}\;\thanks{%
Corresponding author.}}}
\date{}
\maketitle
\begin{abstract}
\msn Semimodules over idempotent semirings like
the max-plus or tropical semiring have much in common with convex cones. 
This analogy is particularly apparent in the case
of subsemimodules of the $n$-fold cartesian product
of the max-plus semiring: it is known that one can separate
a vector from a closed subsemimodule that does not contain it.
We establish here a more general
separation theorem, which applies to any finite collection
of closed subsemimodules with a trivial intersection.
In order to prove this theorem, we investigate
the spectral properties of certain nonlinear operators called here 
idempotent cyclic projectors. These are idempotent analogues 
of the cyclic nearest-point projections known in convex analysis.
The spectrum of idempotent cyclic projectors
is characterized in terms of a suitable extension
of Hilbert's projective metric.
We deduce as a corollary of our
main results the idempotent analogue of Helly's theorem.

\msn {\em Keywords:} Idempotent analysis, tropical semiring,
semimodule, convex geometry, separation, cyclic projections, 
Hilbert's projective metric.

\msn {\em AMS classification (2000):} 52A20 (primary), 06F15, 47H07, 52A01 
(secondary).
\end{abstract}

\section{Introduction}

\msn Some nonlinear problems in optimization theory and mathematical physics
turn out to be linear over semirings with an idempotent
addition $\oplus$ \cite{BCOQ:92}, \cite{C-G:79}, \cite{KM:97<Id>}. 
We recall that the idempotency 
of $\oplus$ means $a\oplus a=a$ for all $a$, and that the role
of this addition is most often played by the operations of taking
maxima or minima.
The search for idempotent
analogues of classical results 
has motivated the development of idempotent
mathematics,
see the recent collection of articles
\cite{LM:05} and also 
\cite{Lit-07<Intro>} for more background.

\msn One of the most studied idempotent semirings 
is the max-plus semiring. It is the 
set $\R\cup\{-\infty\}$ equipped with the operations
of addition $a\oplus b:=\max(a,b)$ and multiplication 
$a\odot b:=a+b$.  
The zero element
$\0$ of this semiring is equal to $-\infty$, and the semiring unity $\1$
is equal to $0$. Some algebraic structures which coincide with the max-plus
semiring (up to isomorphism) have appeared under other names.
In particular, the min-plus or tropical semiring
is obtained by replacing $-\infty$ by $+\infty$ and
$\max(a,b)$ by $\min(a,b)$ above. 
Applying $x\mapsto\exp(x)$ to the max-plus semiring (assume
$\exp(-\infty)=0$), we obtain
the max-times semiring, further denoted by
$\Rmaxm$. It is the set of nonnegative numbers
($\R_+$), equipped with the operations $a\oplus b=\max(a,b)$
and $a\odot b=a\times b$. 
The zero and unit elements of $\Rmaxm$ 
coincide with the usual $0$ and $1$. Our main results 
(Sect.~\ref{s:maxsep}) will be stated over this semiring, as it
makes clearer some analogies with classical convex analysis. 

\msn We shall consider here subsemimodules of the $n$-fold
cartesian product $\cK^n$ of a semiring $\cK$ and,
more generally, of the set $\cK^I$ of functions
from a set $I$ to $\cK$. 
Further examples can be found e.g.\ in \cite{BCOQ:92},
\cite{KM:97<Id>} and \cite{LMSz-01<IFA>}.

\msn In an idempotent semiring, there is a canonical order
relation, for which every element is ``nonnegative''.
Therefore, idempotent semimodules
have much in common with the semimodules over the semiring of nonnegative 
numbers, that is, with {\em convex cones} \cite{Roc:70}. 
One of the first results 
based on this idea is the separation theorem for
convex sets over ``extremal algebras'' proved
by K.~Zimmermann in \cite{KZim-77}. This theorem
implies that a point in $\Rmaxmn$, which does not
belong to a semimodule that is closed in the Euclidean topology, can be
separated from it by an idempotent analogue of a closed halfspace.
Generalizations of this result were obtained in a work by 
S.N.~Samborski{\u{\i}}
and G.B.~Shpiz \cite{SS-92} and in works by G.~Cohen, J.-P.~Quadrat, 
I.~Singer, and
the first author \cite{CGQ-04}, \cite{CGQS-05}. 
In the special case of finitely generated semimodules,
a separation theorem has also been 
obtained by M.~Develin and
B.~Sturmfels in \cite{DS-04},
with a strong emphasis on some combinatorial aspects
of the result.


\msn The main result of this paper, Theorem~\ref{maxsep2}, 
shows that {\em several}
closed semimodules which do not have common nonzero
points can be separated from each other.
This means that for each of these semimodules, we can
select an idempotent halfspace containing it, in such a way
that these halfspaces also do not have common nonzero points.

\msn Even in the case of two semimodules, 
this statement has not been proved in the idempotent literature. Indeed, 
the earlier separation theorems deal with the separation
of a point from an (idempotent) convex set or semimodule,
rather than with the separation of two convex sets or semimodules.
Note that unlike in the classical case, 
separating two convex sets cannot be  
reduced to separating a 
point from a convex set. More precisely, it is easily shown that
two convex sets $A$ and $B$ can be separated if and only if 
the point $0$ can be separated from their Minkowski difference $A-B$, in 
classical convex geometry. In idempotent geometry, an analogue
of Minkowski difference can still be defined, consider
$A\ominus B=\{x\mid\exists\, b\in B:\, x\oplus b\in A\}$.  
However, due to the idempotency of the addition, we cannot recover
a halfspace separating $A$ and $B$ from a halfspace separating $0$
from $A\ominus B$.

\msn In order to prove the main result, Theorem~\ref{maxsep2}, 
we investigate the spectral
properties of idempotent cyclic
projectors. By idempotent cyclic projectors
we mean finite compositions of certain nonlinear projectors on
idempotent semimodules.
The continuity and homogeneity of these nonlinear projectors
enables us to apply to their compositions, i.e.\ 
to the cyclic projectors, some results from non-linear
Perron-Frobenius theory. The main idea is to prove 
the equivalence of the following 
three statements: 1) that the 
semimodules have trivial intersection, 2) that the
separating halfspaces exist, and 3) that the spectral radius
of the associated cyclic projector is strictly 
less than $1$. This
equivalence is established in Theorems~\ref{gensep} and \ref{maxsep},
which deal with the special case of archimedean semimodules, i.e.\ 
semimodules containing at least one positive vector.
As an ingredient of the proof, we use a nonlinear 
extension of Collatz-Wielandt's 
theorem obtained by R.D.~Nussbaum~\cite{Nus-86}.
To derive the main
separation result, Theorem~\ref{maxsep2}, we show that for any
collection of trivially intersecting semimodules, there is a collection of
trivially intersecting {\em archimedean} semimodules, such that every
semimodule from the first collection is contained in an archimedean 
semimodule from the second collection.

\msn We also show in Theorems~\ref{maximum}
and \ref{dhincr} that the orbit of an eigenvector of
a cyclic projector maximizes a certain
objective function. We call this maximum
the Hilbert value of semimodules, as it is a natural
generalization of Hilbert's projective metric, and 
characterize the spectrum of cyclic projectors 
in terms of these Hilbert values
(Theorem \ref{C-G}).

\msn The projectors on idempotent semimodules, which constitute
the cyclic projectors considered here, have been studied
by R.A. Cuninghame-Green, see \cite{C-G-76} and
\cite{C-G:79}, Chapter 8, where they appear as $AA^*$-products.  
The geometrical properties of these projectors
have been used in \cite{CGQ-04,CGQS-05} to establish separation theorems. 
The same operators
have also been studied by G.L.~Litvinov, V.P.~Maslov and 
G.B.~Shpiz, who obtained
in \cite{LMSz-01<IFA>} idempotent analogues of several
results from functional analysis, including
the analytic form of the Hahn-Banach theorem.

\msn The idempotent cyclic projectors
have been introduced, in the case of two semimodules, by 
P.~Butkovi\v{c} and R.A.~Cuninghame-Green in \cite{C-GB-03}, where 
these operators
give rise to an efficient (pseudo-polynomial)
algorithm for finding a point in 
the intersection of two finitely
generated subsemimodules of $\Rmaxmn$. 
In convex
analysis and optimization theory, an analogous role is played by
the cyclic nearest-point projections on convex sets \cite{BBL-97}.

\msn As a corollary of Theorems \ref{maxsep} and \ref{maxsep2}, we
deduce a max-plus analogue
of Helly's theorem. 
This result has also been obtained, 
with a different proof, by F.~Meunier and the first author \cite{Meu-06}.

\msn Our main results apply to subsemimodules of $\Rmaxmn$. 
Some of our results still hold in a 
more general setting, see Sect.~\ref{s:gensep}. However, the separation of 
several semimodules in such a generality remains an open question.

\msn The results of this paper are presented as follows. 
Sect.\ \ref{s:pointsep} describes the main assumptions that are
satisfied by the semimodules of the paper. Besides that, it is occupied
by some basic notions and facts that will be used further.
Sect.\ \ref{s:gensep} is devoted to the results obtained in the
most general setting, with respect to the assumptions of 
Sect.\ \ref{s:pointsep}. The main results for the case $\Rmaxmn$ are obtained 
in Sect.\ \ref{s:maxsep}. These results include separation of several 
semimodules and characterization of the spectrum of cyclic projectors.

\section{Preliminary results on projectors and separation
}
\label{s:pointsep}

\msn  We start this section with some details
concerning the role of partial order in idempotent algebraic structures. 
For more background, we refer the reader to e.g.~\cite{BCOQ:92,C-G:79,
LMSz-01<IFA>}. 

\msn  The idempotent addition $\oplus$ defines the
canonical order relation $\leqoplus$ on the semiring $\cK$
by the rule
$\lambda\oplus \mu=\mu\Leftrightarrow \lambda\leqoplus \mu$ for 
$\lambda,\mu\in\cK$. The idempotent sum $\lambda\oplus\mu$ is equal to
the least upper bound $\sup(\lambda,\mu)$
with respect to the order $\leqoplus$. 
The idempotent sum of an arbitrary subset
is defined to be the least upper bound of this subset, if this least upper
bound exists. In a semimodule $\cV$,
we define the order relation $\leqoplusV$ in the same way. 
The relation $\lambda\leqoplus\mu$ between $\lambda,\mu\in\cK$
implies $\lambda x\leqoplusV\mu x$ for all $x\in\cV$. When $\cV=\cK^n$ and
$\cK=\Rmaxm$, the order 
$\leqoplus$ coincides with the usual 
linear order on $\Rp$, and the order $\leqoplusV$ coincides with
the standard pointwise order on $\Rn$. For this reason, 
we will write
$\leq$ instead of $\leqoplus$ and $\leqoplusV$, in the sequel.

\msn A semiring or a semimodule will be called {\em $b$-complete} 
(see \cite{LMSz-01<IFA>}), if it is
closed under the sum (i.e.\ the supremum) of any subset
bounded from above, and the multiplication  
distributes over such sums.  If the least upper bound $\oplus$ exists 
for all subsets bounded
from above, then the greatest lower bound $\wedge$ exists for all subsets
bounded from below. Consequently, the greatest
lower bound exists for any subset of a $b$-complete
semiring or a semimodule, since such a subset
is bounded from below by $\0$. 

\msn Also note that if $\cK$ is a $b$-complete semiring, and 
the set $\cK\setminus\{\0\}$ is a multiplicative group, 
then this group is abelian by Iwasawa's theorem \cite{Bir:67}.
A semiring $\cK$ such that the set $\cK\setminus\{\0\}$ is an abelian 
multiplicative group is called an {\em idempotent semifield}.

\msn  We shall consider
semirings $\cK$
and semimodules $\cV$ over $\cK$ 
that satisfy the following assumptions:

\msn $(A0)$: the semiring $\cK$ is a $b$-complete idempotent 
semifield, 
and the semimodule $\cV$ is a $b$-complete semimodule
over $\cK$;

\msn $(A1)$: for all elements $x$ and $y\ne\0$ from $\cV$,
the set $\{\lambda\in\cK\mid \lambda y\leq x\}$ is bounded from above.



\msn Assumptions $(A0,A1)$ imply that the operation
\begin{equation}
\label{def:div}
x/y=\max\{\lambda\in\cK\mid \lambda y\leq x\}.
\end{equation}
is defined for all elements $x$ and $y\ne\0$ from $\cV$.
The following can be viewed as another definition of the operation
$/$ equivalent to (\ref{def:div}):
\begin{equation}
\label{resid}
\lambda y\leq x\Leftrightarrow \lambda\leq x/ y.
\end{equation}
In the case $\cV=\cK^I$,
\begin{equation}
\label{div:simple}
x/ y=\bigwedge_{i\colon y_i\ne\0} x_i/ y_i.
\end{equation}

\msn The operation $/$ has the following properties:
\begin{equation}
\label{div:additive}
(\bigwedge_{\alpha} x_{\alpha})/ y=
\bigwedge_{\alpha} (x_{\alpha}/ y),\quad 
(x/\bigoplus_{\alpha} y_{\alpha})=
\bigwedge_{\alpha}(x/ y_{\alpha})
\end{equation}
\begin{equation}
\label{div:hom}
(\lambda x)/ y=\lambda(x/ y)\ \forall\lambda,\quad
y/(\lambda x)=\lambda^{-1}(y/ x)\ \forall\lambda\ne\0.
\end{equation}

\msn We also need the following lemma.
\begin{lemma}
\label{x/x=1}
Under $(A0,A1)$, $x/x=\1$ for all nonzero vectors $x\in\cV$. If
$\lambda x=x$ for a nonzero vector $x\in\cV$, then $\lambda=\1$.
\end{lemma}
\begin{proof}
The inequality $x\leq x$ implies that $x/x\geq\1$, see (\ref{def:div}).
On the other hand, we have that $(x/x)x\leq x$. Multiplying this by
$x/x$, we obtain that $(x/x)^2 x\leq(x/x)x\leq x$, hence $(x/x)^2\leq x/x$
and $x/x\leq\1$. Thus $x/x=\1$.

\msn If $\lambda x=x$ for some $x\ne\0$, then 
$\lambda(x/x)=(\lambda x)/x=x/x$ and so $\lambda=\1$.
\end{proof}\eproof

\begin{definition}
\label{d:subsem}
A subsemimodule $V$ of $\cV$ is a {\em $b$-(sub)semimodule}, 
if $V$ is closed under the sum of any of its subsets bounded
from above in $\cV$.
\end{definition}

\msn Let $V$ be a $b$-subsemimodule of
the semimodule $\cV$. Consider the operator $P_V$ defined by
\begin{equation}
\label{projector}
P_V(x)=\max\{u\in V\mid u\leq x\},
\end{equation}
for every element $x\in\cV$. Here we use ``$\max$'' to
indicate that the least upper bound belongs to the set. The operator $P_V$
is a {\em projector} onto the subsemimodule $V$, as $P_V(x)\in V$ for
any $x\in\cV$ and $P_V(v)\in V$ for any $v\in V$. In principle,
$P_V$
can be defined for all subsets of $\cV$, if we write $\sup$ instead
of $\max$ in (\ref{projector}), but 
then $P_V$ may not be a projector on $V$.

\begin{definition}
\label{d:elem}
A subsemimodule $V$ of $\cV$ is called {\em elementary},
if $V=\{\lambda y\mid \lambda\in\cK\}$ for some $y\in\cV$.
The projector onto such a semimodule is also called {\em elementary}.
\end{definition}

\msn Assumptions $(A0,A1)$ imply that elementary semimodules are
$b$-semimodules.  For the elementary semimodule
$V=\{\lambda y\mid \lambda\in\cK\}$, the projector $P_V$ 
is given by $P_V(x)=(x/ y)y$,
and this fact can be generalized as follows.  

\begin{proposition}
\label{ortproj2}
If $V$ is a $b$-subsemimodule of $\cV$ and
$P_V(x)=\lambda y$ for some $\lambda\in\cK$ and $x,y\in\cV$,
then
$P_V(x)=(x/ y)y$.
\end{proposition}
\begin{proof} If $V$ is a $b$-semimodule, 
then $y\in V$, and $(x/ y) y\leq x$  
implies that $(x/ y)y\leq P_V(x)=\lambda y$. On the other hand, 
$\lambda y\leq x$ implies $\lambda\leq x/ y$
and $\lambda y\leq (x/ y)y$.
\end{proof}\eproof
Note that $P_V$
is isotone with respect to inclusion:
\begin{equation}
\label{incl-iso}
U\subset V\Rightarrow P_U(x)\leq P_V(x)\ \text{for all $x$}.
\end{equation}
It is also homogeneous and isotone:
\begin{equation}
\label{e:projprops}
P_V(\lambda x)=\lambda P_V(x),\quad x\leq y\Rightarrow P_V(x)\leq P_V(y).
\end{equation}
We remark that the operator 
$P_V$ is in general {\em not} linear with respect
to $\oplus$ or $\wedge$ operations, even in the case $\cV=\Rmaxmn$.

\msn In idempotent geometry, the role of halfspace is played
by the following object.

\begin{definition}
\label{d:halfs}
A set $H$ given by
\begin{equation}
\label{halfs0}
H=\{x\mid u/x \geq v/x\}\cup\{\0\}
\end{equation}
with $u,v\in \Rmaxmn$, $u\leq v$, 
will be called {\em (idempotent) halfspace}.
\end{definition}

\msn Properties (\ref{div:additive}) and (\ref{div:hom}) of the operation 
$/$ imply that any halfspace is a semimodule. If $\cV=\cK^I$,
then we can use 
(\ref{div:simple}) and then
\begin{equation}
\label{h:simple}
H=\{x\mid\bigwedge_{i\colon x_i\ne\0} u_i x_i^{-1}\geq
\bigwedge_{i\colon x_i\ne\0} v_i x_i^{-1}\}\cup\{\0\}.
\end{equation}
If $\cV=\cK^n$ and all coordinates of $u$ and $v$ are nonzero,  
then we have that
\begin{equation}
\label{reghalfspace}
H=\{x\mid\bigoplus_{\ldotsn} 
x_i u_i^{-1}\leq
\bigoplus_{\ldotsn} x_i v_i^{-1}\}.
\end{equation}
Such idempotent halfspaces formally resemble the
closed homogeneous halfspaces of the 
finite-dimensional convex geometry \cite{Roc:70}.

\msn Since the operation $/$ is 
isotone with respect to the first argument, we can replace the
inequalities in (\ref{d:halfs}), (\ref{h:simple})
and (\ref{reghalfspace}) by the equalities. For instance, 
definition (\ref{d:halfs}) can be rewritten as
\begin{equation}
\label{halfspace1}
H=\{x\mid u/ x= v/ x\}\cup\{\0\},
\end{equation}
where $u\leq v$.

\msn The present paper is concerned with 
the separation of several 
$b$-semimodules, whereas the 
separation theorems which have 
been established previously,
like the ones of \cite{CGQ-04,CGQS-05}, 
deal with the separation of one point 
from a semimodule. 
For the convenience of the reader, 
we next state a theorem, which is a variant of a 
separation theorem of \cite{CGQ-04}. The difference is in that 
we deal with 
$b$-complete semimodules
rather than with complete semimodules. 
Both results are closely related with
the idempotent Hahn-Banach theorem of \cite{LMSz-01<IFA>}.
 

\begin{theorem}
\label{separation}
{\rm (Compare with \cite{CGQ-04}, Theorem~8)}
Let $V$ be a $b$-subsemimodule of $\cV$ and let $u\notin V$. 
Then the halfspace
\begin{equation}
\label{halfspace}
H=\{x\mid P_V(u)/x\geq u/ x\}\cup\{\0\}
\end{equation}
contains $V$ but not $u$.
\end{theorem}
\begin{proof} Take a nonzero vector $x\in V$ 
(the case $x=\textbf 0$ is trivial). Since $(u/x)x\leq u$, we
have $(u/x)x\leq P_V(u)$, which is by (\ref{resid}) equivalent to
$u/x\leq P_V(u)/x$. Hence $V\subseteq H$.

\msn Take $x=u$ and assume that $P_V(u)/u\geq u/u=\1$. 
This is equivalent to 
$u\leq P_V(u)$ and hence to
$u=P_V(u)$. Since $V$ is a $b$-semimodule, we have that $u\in V$, 
which is a contradiction. Hence $u\notin H$.
\end{proof}\eproof

\begin{definition}
\label{def:preceq}
Consider the preorder relation $\preceq$
defined by
\begin{equation}
\label{preceq}
x\preceq y\Leftrightarrow y/ x>\0.
\end{equation}
We say that $x$ and $y$ are
{\em comparable}, 
and we write $x\sim y$,
if $x\preceq y$ and 
$y\preceq x$.  Equivalently,
\begin{equation}
\label{sim}
x\sim y\Leftrightarrow (x/y)(y/x)>\0.
\end{equation}
\end{definition}

\msn Note that if $y=\lambda x$ with $\lambda\ne\0$, then
$y\sim x$, and that the inequality $x\leq y$, if $x\ne\0$, 
implies that $x\preceq y$.
In particular, $P_V(x)\preceq x$ for any nonzero $x\in\cV$ 
and any semimodule $V$, provided that $P_V(x)$ is nonzero. 

\msn When $\cV=\cK^n$, 
comparability can be characterized
in terms of supports. 
Recall that the {\em support}
of a vector $x$ in $\cK^n$
is defined by 
$\supp(x)=\{i\mid x_i\ne 0\}$.
It can be checked that for all 
$x,y\in\cK^n$, we have $x\preceq y$ 
iff $\supp(x)\subset \supp(y)$, 
and so, $x\sim y$ 
iff $\supp(x)=\supp(y)$.

\begin{proposition}
\label{simple1}
Let $x\in\cV$ be a nonzero vector and let $V\subseteq\cV$ be a
$b$-semimodule containing a nonzero vector $y$.
If $y\preceq x$, then
$P_V(x)$ is nonzero, and $y\preceq P_V(x)\preceq x$. If $y\sim x$,
then $P_V(x)\sim x$.
\end{proposition}
\begin{proof} 
By the definition of $/$ and by (\ref{preceq}),
there exists $\alpha$ such that $\alpha y\leq x$.
Then $\alpha y\leq P_V(x)$, hence $P_V(x)$ is nonzero and  
$y\preceq P_V(x)$. 
\end{proof}\eproof

\begin{proposition}
\label{l<m}
Let $F$ be an isotone and homogeneous operator, let $\lambda,\mu$ be
arbitrary scalars from $\cK$
and let $v$ and $u$ be nonzero vectors such that $v\prec u$.
Suppose that one of the following is true:
\begin{itemize}
\item[1.] $Fv\geq\mu v$ and $Fu=\lambda u$;
\item[2.] $Fv=\mu v$ and $Fu\leq\lambda u$.
\end{itemize}
Then $\mu\leq\lambda$.
\end{proposition}
\begin{proof}  Applying $F$ to the inequality
$(u/ v) v\leq u$ and using any of the given conditions, 
we obtain that
$(u/ v) \mu v\leq \lambda u$. If $\lambda=\0$,
then $\mu=\textbf{0}$. If $\lambda$ is invertible, then by (\ref{resid})
$(u/ v)\mu\lambda^{-1}\leq u/ v$. Cancelling $u/ v$,
we get $\mu\leq\lambda$.
\end{proof}\eproof

\msn  Properties (\ref{div:additive}) and (\ref{div:hom}) 
imply that the sets
 $\{x\mid x\preceq y\}$, $\{x\mid x\succeq y\}$
and hence $\{x\mid x\sim y\}$ are subsemimodules of $\cV$. For
any semimodule $V\subset\cV$ and any vector $y\in\cV$, we
define 
\begin{equation}
\label{vy}
V^y=\{x\in V\mid x\preceq y\},
\end{equation}
which is a subsemimodule of $V$. When $\cV=\cK^n$,
$V^y$ is uniquely determined by the support $M$ of $y$. For this
reason, for all $M\subseteq \{1,\ldots,n\}$, we set
\begin{equation}
\label{vm}
V^M=\{x\in V\mid \text{supp}(x)\subset M\}.
\end{equation}

\begin{definition}
\label{d:arch}
A vector $x\in\cV$ is called 
{\em archimedean}, 
if $y\preceq x$ for all $y\in\cV$.
A semimodule $V\subseteq\cV$ 
is called
{\em archimedean}, if it contains an archimedean vector.
A halfspace will be called
{\em archimedean},
if both vectors defining it (e.g.\ $u$ and $v$ in 
(\ref{halfs0})) are archimedean.
\end{definition}

\msn Of course, Def.\ \ref{d:arch} makes sense only in the case
when $\cV$ satisfies the following assumption:

\msn $(A2)$: the semimodule $\cV$ has an archimedean vector.

\msn This assumption is satisfied by the semimodules $\cV=\cK^n$ 
(we are also assuming $(A0,A1)$).
In this case, 
archimedean halfspaces have been written explicitly in
(\ref{reghalfspace}).

\section{Cyclic projectors and separation theorems: general results}
\label{s:gensep}
\msn In this section we study cyclic projectors, that is,
compositions of projectors
\begin{equation}
\label{comp-proj}
P_{V_k}\cdots P_{V_1},
\end{equation}
where $V_1,\ldots,V_k$ are $b$-subsemimodules of $\cV$. We assume $(A0,A1)$,
which means in particular that $\cK$ is an idempotent semifield, and
state general results concerning cyclic projectors and separation properties.
For the notational convenience, 
we will write $P_t$ instead of
$P_{V_t}$. We will also adopt a convention of cyclic numbering
of indices of projectors
and semimodules, so that
$P_{l+k}=P_l$ and $V_{l+k}=V_l$ for all $l$.

\begin{definition}
\label{d:hvalue}
Let $x^1,\ldots,x^k$ be nonzero elements of $\cV$. 
The value
\begin{equation}
\label{e:hvalue}
d_{\Hilb}(x^1,\ldots,x^k)=(x^1\odiv x^2)\ (x^2\odiv x^3)\ldots
(x^k\odiv x^1).
\end{equation}
will be called {\em the Hilbert value} of $x^1,\ldots,x^k$.
\end{definition}

\msn It follows from Def.~\ref{def:preceq} that 
$d_{\Hilb}(x^1,\ldots,x^k)\ne\0$ if and only if
all vectors $x^1,\ldots,x^k$ are comparable. One can show that
$d_{\Hilb}(x^1,\ldots,x^k)\leq\1$. This inequality is an equality if and
only if $x^1,\ldots,x^k$ differ from each other only by scalar multiples.
The Hilbert value is invariant under multiplication of any of its
arguments by an invertible scalar, and by cyclic
permutation of its arguments.

\msn The Hilbert value of
two vectors $x^1,x^2$
was studied in \cite{CGQ-04}. For two comparable
vectors in $\Rmaxmn$, that is, for two vectors with common
support $M$ it is given by
\begin{equation}
\label{e:hvalue3}
d_{\Hilb}(x^1,x^2)=\min_{i,j\in M}
(x^1_i(x^2_i)^{-1}x^2_j (x^1_j)^{-1}),
\end{equation}
so that $-\log(d_{\Hilb}(x^1,x^2))$ coincides with 
Hilbert's projective metric
\begin{equation}
\label{e:hdist}
\delta_{\Hilb}(x^1,x^2)=\log(\max_{i,j\in M}
(x^1_i(x^2_i)^{-1}x^2_j (x^1_j)^{-1}))=-\log(d_{\Hilb}(x^1,x^2)).
\end{equation}

\begin{definition}
\label{d:hvalue-sem}
The {\em Hilbert value} of $k$ subsemimodules 
$V_1,\ldots,V_k$ of $\cV$ is defined by 
\begin{equation}
\label{hvalue-sem}
d_{\Hilb}(V_1,\ldots,V_k)=\sup_{x^1\in V_1,\ldots,x^k\in V_k}
d_{\Hilb}(x^1,\ldots,x^k)
\end{equation}
\end{definition}

\begin{theorem}
\label{maximum}
Suppose that
the operator $P_k\cdots P_1$ has an 
eigenvector $y$ with eigenvalue $\lambda$. Then
\begin{equation}
\label{strongmax}
\lambda=\max_{x^1\in V_1^y,\ldots,x^k\in V_k^y} d_{\Hilb}(x^1,\ldots,x^k)=
d_{\Hilb}(\Bar{x}^1,\ldots,\Bar{x}^k),
\end{equation}
where $\Bar{x}^i=P_i\cdots P_1 y$.
\end{theorem}
\begin{proof} Note that $\Bar{x}^i$, for any $i$, 
is an eigenvector of $P_{i+k}\cdots P_{i+1}$ and that all
these vectors are comparable with $y$.
Further, let $x^1,\ldots, x^k$ be arbitrary elements of
$V^y_1,\ldots,V^y_k$, respectively, and let  
$\alpha_1,\ldots,\alpha_k$ be scalars such that
\begin{equation}
\label{alphas}
\begin{array}{lcl}
\alpha_1 x^2&\leq &P_2 x^1,\\
&\vdots\\
\alpha_{k-1} x^k&\leq& P_k x^{k-1},\\
\alpha_k x^1&\leq& P_1 x^k,\\
\end{array}
\end{equation}
Take the last inequality. Applying $P_2$
to both sides and using the first inequality, we have that 
$\alpha_1\alpha_k x^2\leq P_2 P_1 x^k$. Further, we apply $P_3$
to this inequality and use the inequality $\alpha_2 x^3\leq P_3 x^2$.
Proceeding in the same manner, we finally obtain 
\begin{equation}
\label{alphas2}
\alpha_1\ldots\alpha_k x^k\leq P_k\cdots P_1 x^k.
\end{equation}
It follows from Prop.~\ref{l<m} that $\alpha_1\cdots\alpha_k\leq\lambda$.
We take
$\alpha_i=x^i/ x^{i+1}$ for $i=1,\ldots, k-1$,
and $\alpha_k=x^k\odiv x^1$. This leads to
\begin{equation}
\label{weakmax}
d_{\Hilb}(V_1^y,\ldots,V^k_y)\leq\lambda.
\end{equation}
Note that this inequality is true if $V_1,\ldots,V_k$ are not 
$b$-semimodules.
Applying Prop.~\ref{ortproj2} we have that 
$\lambda y=d_{\Hilb}(\Bar{x}^1,\ldots,\Bar{x}^k) y$. By Lemma
\ref{x/x=1} we can cancel $y$, and the observation that $\Bar{x}^i\in V_i^y$
for all $i$ yields the desired equality.
\end{proof}\eproof
 
\msn The situation when $P_k\cdots P_1$ has an eigenvector with 
nonzero eigenvalue
occurs, if at least one of the semimodules $V_1,\ldots,V_k$ is
elementary, that is, generated by a single vector 
$x^i$, and if all other semimodules have vectors comparable with $x^i$.
In this case $P_k\cdots P_{i+1}x^i$ is the
only eigenvector of $P_k\cdots P_1$ with nonzero
eigenvalue. 

\msn To obtain the following lemma, we use Prop. \ref{simple1}.
\begin{lemma}
\label{simple3}
Let $x^1\in V^1$ and $x^i=P_i x^{i-1}$ for $i=2,\ldots,k$. Then, the
Hilbert value
$d_{\Hilb}(x^1,\ldots,x^k)$ is not equal to $\0$ if and only if 
$V_2,\ldots,V_k$ have vectors comparable with $x^1$.
\end{lemma}

\begin{theorem}
\label{dhincr}
Suppose that the vectors
$x^i,\ i=1,\ldots$ are such that $x^1\in V_1$ and
$x^i=P_i x^{i-1}$ for $i=2,\ldots$. Then 
$d_{\Hilb}(x^{l+1},\ldots,x^{l+k})$ is nondecreasing with $l$
so that the following inequalities hold for all $l$:
\begin{equation}
\label{dhineqs}
d_{\Hilb}(x^1,\ldots,x^k)\leq
d_{\Hilb}(x^{l+1},\ldots,x^{l+k})\leq\1.
\end{equation}
\end{theorem}
\begin{proof}
As $V_i$ are $b$-semimodules, $x^i\in V_i$ for all $i$. 
If the Hilbert value is $\0$ for all $l$, then
there is nothing to prove. So we assume that there exists
a least $l=l_{\min}$ for which the Hilbert value 
$d_{\Hilb}(x^{l},\ldots, x^{l+k-1})$ is nonzero. As 
it is nonzero, by Lemma \ref{simple3}, all
$x^l,\ldots, x^{l+k-1}$ are comparable. By
Prop. \ref{simple1}, 
$x^{l+k}$ is also comparable with them, and the same is true about the rest of
the sequence, hence $d_{\Hilb}(x^l,\ldots, x^{l+k-1})$ is nonzero for all
$l\geq l_{\min}$. Now we take any $l\geq l_{\min}$ 
and consider the composition 
\begin{equation}
\label{newcomp}
P_{l+k} P'_{l+k-1}\cdots P'_{l+1},
\end{equation}
where $P'_i$, for
$i=l+1,\ldots,l+k-1$, are elementary projectors onto the
semimodules generated by
$x^i$. The operator (\ref{newcomp}) has an eigenvector $x^{l+k}$.
By Theorem \ref{maximum}
\begin{equation}
\label{use2}
d_{\Hilb}(x^l,\ldots,x^{l+k-1})\leq\\[1.5 ex]
\max_{y\in V_l,\ y\preceq x^{l+k}} 
d_{\Hilb}(x^{l+1},\ldots,x^{l+k-1},y)=d_H(x^{l+1},\ldots,x^{l+k}).
\end{equation}
for all $l=1,\ldots$. 
\end{proof}\eproof

The following theorem assumes the existence of archimedean vectors $(A2)$.
\begin{theorem}
\label{gensep}
Suppose that $P_k\cdots P_1$ has an archimedean eigenvector $y$
with nonzero eigenvalue $\lambda$. The following are equivalent:
\begin{itemize}
\item[1.] there exist an archimedean vector $x$ and a scalar $\mu<\1$
such that
$P_k\cdots P_1 x\leq\mu x$;
\item[2.] for all $i=1,\ldots,k$ there exist archimedean halfspaces
$H_i$ such that $V_i\subseteq H_i$ and $H_1\cap\cdots\cap H_k=\{\0\}$;
\item[3.] $V_1\cap\cdots\cap V_k=\{\0\}$;
\item[4.] $\lambda<\1$.
\end{itemize}
\end{theorem}

\begin{proof} $1\Rightarrow 2$:
Denote $x^0=x$ and 
$x^i=P_i\cdots P_1 x^0$. Note that all the $x^i$
are also archimedean by Prop. \ref{simple1}. For all $i=1,\ldots,k$
we have that
\begin{equation}
\label{vinh}
V_i\subseteq\{u\colon x^{i-1}\odiv u=x^i\odiv u\}=H_i.
\end{equation}
Indeed, if $x^{i-1}=x^i$, then $H_i$ coincides
with the whole $\cV$.  If $x^{i-1}\ne x^i$, which means that
$x^i\notin V_{i-1}$,
then the inclusion in (\ref{vinh}) follows from Theorem \ref{separation}.
Assume that there exists a nonzero vector $u$ which belongs to every $H_i$.
Then $x^k\odiv u=x\odiv u$. But
$x^k\odiv u\leq(\mu x)\odiv u\leq x\odiv u$, 
hence $\mu(x\odiv u)=(\mu x)\odiv u=x\odiv u$. Cancelling
$x\odiv u$, we get $\mu=\1$ which contradicts
$1.$ The implication is proved.

\msn $2\Rightarrow 3$: Immediate.

\msn $3\Rightarrow 4$: By the conditions of this theorem,
$P_k\cdots P_1$ has an eigenvector $y$ with eigenvalue $\lambda$.
As any vector is greater than or equal to its image by the projector $P_i$, 
we have that
$\lambda\leq\1$.  Assume that $\lambda=\1$.
Then the inequalities
\begin{equation}
\label{projineqs2}
P_k\cdots P_1 y\leq P_{k-1}\cdots P_1 y\leq\ldots\leq y
\end{equation}
turn into equalities, and $y$ is a common vector of
$V_1,\ldots,V_k$, which contradicts $3.$

\msn $4\Rightarrow 1$: Take $x=y$. 
\end{proof}\eproof

\begin{figure}[htpb]
\begin{center}
\begin{tabular}{ccc}
\input{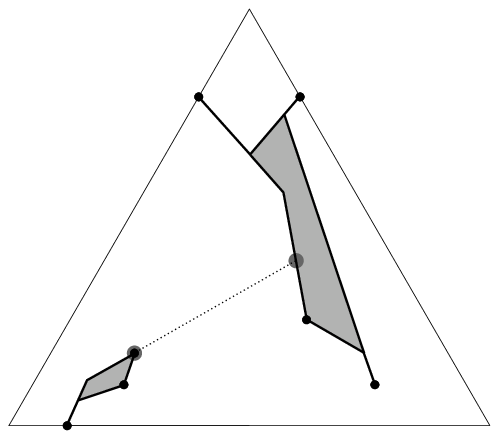}&&\input{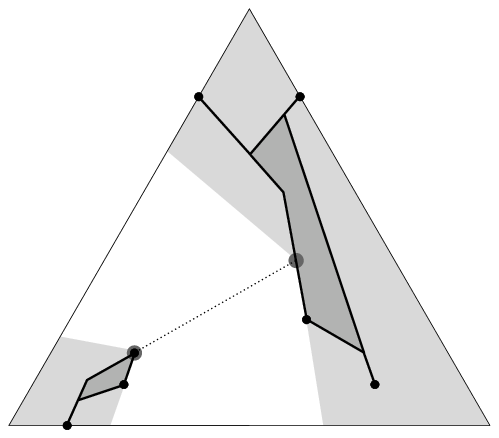}
\end{tabular}
\end{center}
\caption{Two semimodules (left) with separating halfspaces (right)}
\label{fig1}
\end{figure}

\msn To illustrate Theorem~\ref{gensep}, consider the matrices
\begin{equation}
\label{abmatr}
A= \left(\begin{array}{cccc} 0&0&0&-\infty\\ 1&2&-\infty&1\\ 0&-1&2&3
\end{array}\right),
\qquad 
B= \left(\begin{array}{ccc} 3&2&2\\ 0&0&0\\ -\infty&0&-1
\end{array}\right) \enspace .
\end{equation}
Let $a^i$ and $b^i$ denote the $i$-th
column of $A$ and $B$, respectively. 
For all vectors $x=(x_1,\ldots,x_n)$ and 
$\beta>0$, we denote by $\exp(\beta x)$ the
vector of the same size with entries $\exp(\beta x_j)$. 
We define $V_1$ (resp.\ $V_2$) to be the subsemimodule of $\Rmaxm^3$
generated by the vectors $\exp(\beta a^i)$ for $1\leq i\leq 4$ 
(resp.\ $\exp(\beta b^i)$ for $1\leq i\leq 3$).
The discussions which follow are independent of the choice
of the scaling parameter $\beta>0$, which is adjusted to
make Figure~\ref{fig1} readable (we took $\beta=2/3$).
The two semimodules $V_1$, $V_2$ and their generators are represented as
follows at the left of the figure.
Here, a non-zero vector $w=(w_1,w_2,w_3)\in\Rmaxm^3$ is represented
by the point of the two dimensional simplex
which is the barycenter with weights $w_j$ of the three vertices
of this simplex. The generators $a^i$ and $b^i$ correspond
to the bold dots. The semimodules $V_1$ and $V_2$ 
correspond to the two medium grey regions, together with the
bold broken segments joining the generators to each of these regions.

\msn Since the entries of $x^0:=b^2=\exp(\beta(2,0,0))\in V_2$ 
are nonzero, the vector $x^0$
is archimedean, and one can check, using the explicit
formula of the projector (Theorem~5 of~\cite{CGQ-04}),
that $x^0$ is an eigenvector
of $P_2P_1$. Indeed,
\begin{equation}
\label{x1=p1x0}
x^1:=P_1x^0= \exp(\beta(-1,0,0)^T)
\end{equation}
and
\begin{equation}
\label{x2=p2x1}
x^2:=P_2x^1=\exp(\beta(-1,-3,-3)^T)=\exp(-3 \beta) x^0  \enspace .
\end{equation}
The halfspaces constructed in the proof of Theorem~\ref{gensep} 
are given by
\begin{equation}
\begin{array}{l@{{}={}}l}
H_1&\{u\mid x^0/u=x^1/u\}\\
   &\{u\mid \min(\exp(2\beta)/u_1,1/u_2, 1/u_3)=\min(\exp(-\beta)/u_1,1/u_2,1/u_3)\}\\
   &\{u\mid \max(u_2,u_3)\geq \exp(\beta)u_1\}
\end{array}
\end{equation}
and
\begin{equation}
\begin{array}{l@{{}={}}l}
H_2&\{u\mid x^1/u=x^2/u\}\\
   &\{u\mid \min(\exp(-\beta)/u_1,1/u_2,1/u_3)= \min(\exp(-\beta)/u_1,\exp(-3\beta)/u_2, \exp(-3\beta)/u_3)\}\\
   &\{u\mid u_1\geq \exp(2\beta)\max(u_2,u_3)\}.
\end{array}
\end{equation}
These two halfspaces are represented by the zones in light gray
(right). The proof of Theorem~\ref{gensep} shows that their intersection
is zero (meaning that it is reduced to the zero vector).

\section{Cyclic projectors and separation theorems in $\Rmaxmn$}
\label{s:maxsep}
\msn In $\Rmaxmn$, it is natural to consider semimodules that are closed
in the Euclidean topology. One can easily show that such semimodules are 
$b$-semimodules.
Theorem 3.11 of \cite{CGQS-05} implies that projectors onto
closed subsemimodules of $\Rmaxmn$ are continuous.

\msn In order to relax the assumption concerning archimedean
vectors in Theorem \ref{gensep}, we shall use some results from
nonlinear spectral theory, that we next recall. 
By Brouwer's fixed point theorem, a continuous
homogeneous operator $x\mapsto Fx$ that maps $\Rpn$ to itself
has a nonzero eigenvector. This allows us to define
the nonlinear spectral radius of $F$, 
\begin{equation}
\label{specrad}
\rho(F)=\max\{\lambda\in \R_+\mid \exists x\in 
(\R_+^n)\setminus 0,\; Fx=\lambda x\} \enspace .
\end{equation}
Suppose in addition that $F$ is isotone.
Then it can be checked 
that if $Fx=\lambda x$, $Fy=\mu y$, 
and if $x$ and $y$ are comparable, then $\lambda=\mu$. It follows that the number of eigenvalues
of $F$ is bounded by the number of nonempty supports of vectors
of $\R_+^n$, i.e.\ by $2^n-1$. This implies in particular that the maximum
is attained in (\ref{specrad}). We shall need the following nonlinear
generalization of the Collatz-Wielandt formula for the spectral
radius of a nonnegative matrix.

\begin{theorem}
\label{t:Nussbaum}
{\rm (R.D.~Nussbaum, Theorem~3.1 of \cite{Nus-86})}
For any isotone, homogeneous, and continuous map $F$ from
$\Rpn$ to itself, we have:
\begin{equation}
\label{Nussbaum}
\rho(F)=\inf_{x\in(\Rp\backslash\{0\})^n}\max_{1\leq i\leq n} [F(x)]_i x_i^{-1}.
\end{equation}
\end{theorem}
This result implies that 
the spectral radius of such operators is isotone:
if $F(x)\leq G(x)$ for any $x\in\Rpn$, then 
$\rho(F)\leq\rho(G)$.

\msn As the projectors on subsemimodules of $\Rmaxmn$ are isotone, 
homogeneous and continuous, so are their compositions, i.e.\ cyclic
projectors. Consequently, we can apply Theorem~\ref{t:Nussbaum} 
to them. 
Hence, if $V'_i,\ i=1,\ldots,k$ and $V_i,\ i=1,\ldots,k$ 
are closed semimodules in $\Rmaxmn$ and such that 
$V'_i\subseteq V_i,\ i=1,\ldots,k$, then 
\begin{equation}
\label{e:monot} 
\rho(P'_k\cdots P'_1)\leq\rho(P_k\cdots P_1),
\end{equation}
as the projectors are isotone with respect to inclusion (\ref{incl-iso}).

\msn In the following theorem we characterize the spectrum of cyclic
projectors in $\Rmaxmn$.

\begin{theorem}
\label{C-G}
Let $V_1,\ldots,V_k$ be closed semimodules
in $\Rpn$. Then the Hilbert value $d_{\Hilb}(V_1,\ldots,V_k)$ is the
spectral radius of $P_k\cdots P_1$.
The spectrum of $P_k\cdots P_1$ is the set of Hilbert values
$d_{\Hilb}(V^M_1,\ldots,V^M_k)$, where $M$ ranges over all
nonempty subsets
of $\{1,\ldots,n\}$. 
\end{theorem}

\begin{proof}  
We first prove that the Hilbert value $d_{\Hilb}(V_1,\ldots, V_k)$
is the spectral radius of the cyclic
projector, and hence an eigenvalue. 
We take $k$ elementary subsemimodules 
spanned by $x^i\in V_i,\ i=1,\ldots,k$ and consider
elementary projectors $P'_i$ onto them. Observe that 
\begin{equation}
\rho(P'_k\cdots P'_1)=d_{\Hilb}(x^1,\ldots,x^k).
\end{equation} 
Denote by $\Bar{x}^0$ an eigenvector of $P_k\cdots P_1$,
associated with the spectral radius, and let 
$\Bar{x}^i=P_i\cdots P_1 \Bar{x}^0$. Then
\begin{equation}
\rho(P_k\cdots P_1)=d_{\Hilb}(\Bar{x}^1,\ldots,\Bar{x}^k).
\end{equation}
By (\ref{e:monot}) it follows that 
$\rho(P_k\cdots P_1)\geq\rho(P'_k\cdots P'_1)$,
that is, 
$d_{\Hilb}(\Bar{x}^1,\ldots,\Bar{x}^k)\geq d_{\Hilb}(x^1,\ldots x^k)$
for any $x^1\in V_1,\ldots, x^k\in V^k$.
Thus, the Hilbert value of $V_1,\ldots V_k$ is the spectral radius of
$P_k\cdots P_1$.

\msn Now consider $d_{\Hilb}(V_1^M,\ldots,V_k^M)$ for arbitrary 
$M\subseteq\{1,\ldots,n\}$. Note that the semimodules $V_1^M,\ldots, V_k^M$
are closed, and denote by $P_1^M,\ldots, P_k^M$ the projectors onto these. 
It is easy to see that $P_i^M(y)=P_i(y)$ for all $i$ and 
all $y$ with $\text{supp}(y)\subseteq M$. It follows that 
$d_{\Hilb}(V_1^M,\ldots,V_k^M)$ is the spectral radius of 
$P_k^M\cdots P_1^M$ and also an eigenvalue of
$P_k\cdots P_1$. 

\msn We have proved that any Hilbert value 
$d_{\Hilb}(V_1^M,\ldots,V_k^M)$ is an eigenvalue of
$P_k\cdots P_1$. The converse statement follows from
Theorem \ref{maximum}.
\end{proof}\eproof

The following three results refine Theorem~\ref{gensep}.
\begin{theorem}
\label{maxsep}
Suppose that $V_1,\ldots,V_k$ are closed archimedean subsemimodules of
$\Rmaxmn$. The following are equivalent:
\begin{itemize}
\item[1.] there exist a positive vector $x$ and a number $\lambda<1$ such
that $P_k\cdots P_1 x\leq\lambda x$;
\item[2.] there exist archimedean halfspaces $H_i$ which contain
$V_i$ and are such that $H_1\cap\cdots\cap H_k=\{\0\}$;
\item[3.] $V_1\cap\cdots\cap V_k=\{\0\}$;
\item[4.] $\rho(P_k\cdots P_1)<\1$.
\end{itemize}
\end{theorem}
\begin{proof} The implications $1\Rightarrow 2$, $2\Rightarrow 3$ и 
$3\Rightarrow 4$ are proved in Theorem \ref{gensep}.
The implication $4\Rightarrow 1$ follows from Equation (\ref{Nussbaum}).
\end{proof}\eproof
\begin{proposition}
\label{include}
Suppose that $V_i,\ i=1,\ldots,k$
are closed semimodules in $\Rmaxmn$ with zero intersection.
Then there exist closed archimedean semimodules
$V'_i,\ i=1,\ldots,k$ with zero intersection and such that
each $V'_i$ contains $V_i$.
\end{proposition}
\begin{proof}
In every semimodule $V_i$, we find a vector $y^i$
with maximal support and such that  
$||y^i||=\max(y^i_1,\ldots,y^i_n)=1$.
For all scalars $\delta>0$, define 
\begin{equation}
\label{e:zi1}
z^i(\delta)=y^i\oplus\delta\bigoplus_{j\notin\text{supp}(y^i)} e^j
\end{equation}
and the semimodules
\begin{equation}
\label{vi1}
V_i(\delta)=\{x\mid x=v\oplus\lambda z^i,\ v\in V'_i,\ \lambda\in\Rp\}.
\end{equation}
These semimodules are closed, as all arithmetical operations are continuous.
We show that for $\delta>0$ small enough 
these semimodules have zero intersection.
Assume by contradiction
that for all $\delta>0$, there exists a nonzero vector 
$u(\delta)$ in the 
intersection 
$V_1(\delta)\cap \cdots \cap V_k(\delta)$. 
After normalizing $u(\delta)$, we may assume that $||u(\delta)||=1$.
For any $i=1,\ldots,k$
and any $\delta$ we have that
\begin{equation}
\label{u1}
u(\delta)=v^i(\delta)\oplus\lambda_i(\delta)y^i\oplus
\lambda_i(\delta)\delta\bigoplus_{j\notin\text{supp}(y^i)} e^j,
\end{equation}
where $v^i(\delta)$ is a vector from $V'_i$
and $\lambda_i(\delta)$ is a scalar.
As $||u(\delta)||=1$ and $||y^i||=1$, we have that 
$\lambda_i(\delta)\leq 1$. So, there exists
a sequence $(\delta_m)_{m\geq 1}$ converging to $0$ such that
for all $1\leq i\leq k$, 
$\lambda_i(\delta_m)$ converges to a limit as $m$ tends to infinity.
Then 
$w:=\lim_{m\to\infty} u(\delta_m)=
\lim_{m\to\infty} v_i(\delta_m)\oplus\lambda_i(\delta_m)y^i$
for all $i$.
As $V_i$ are closed, $w$
belongs to $V_i$ at all $i$. Since $||w||=1$, $w$ is not equal to
$\0$, which is a contradiction.
\end{proof}\eproof

The following is an immediate corollary of Theorem~\ref{maxsep} 
and Proposition~\ref{include}.
\begin{theorem}
\label{maxsep2}
{\rm (Separation theorem)}
If $V_i,\ i=1,\ldots,k$ are closed semimodules
with zero intersection, then 
there exist archimedean halfspaces $H_i,\ i=1,\ldots,k$, which contain
the corresponding semimodules $V_i$ and have zero intersection.
\end{theorem}

The following
separation theorem for two closed
semimodules is a corollary of Theorem \ref{maxsep2}.

\begin{theorem}
\label{2sep}
If $U$ and $V$ are two closed semimodules with zero intersection,
then there exists a closed halfspace $H_U$, which contains
$U$ and has zero intersection with $V$,
and there exists a closed halfspace $H_V$, which contains $V$ and has
zero intersection with $U$.
\end{theorem}

\msn As a consequence of Theorem~\ref{maxsep2}, 
we further deduce a separation
theorem for convex subsets of $\Rmaxmn$.
We recall here some definitions from idempotent convex geometry, 
see e.g. \cite{GK-07<LAA>}
A subset $C\subset \Rmaxmn$ is {\em convex}
if $\lambda u\oplus \mu v\in C$, 
for all $u,v\in C$ and $\lambda,\mu\in \Rmaxm$
such that $\lambda\oplus\mu=1$.

\msn
The {\em recession cone} of a convex set $C$,
$\operatorname{rec}(C)$, is  
the set of vectors $u$ such that $v\oplus \lambda u\in C$
for all 
$\lambda\in \Rmaxm$, 
where $v$ is an arbitrary vector of $C$.
As shown in Prop.~2.6 of~\cite{GK-07<LAA>}, 
if $C$ is closed, 
the recession cone is independent of
the choice of $v$.  
Observe that when $C$ is compact, its recession
cone is zero. 

\msn
A set $H^{\aff}$ given by
\begin{equation}
\label{aff-halfs}
H^{\aff}=\{x\mid u/x \wedge \alpha \geq v/x\wedge \gamma\}
\end{equation}
with $u,v\in \Rmaxmn$, $u\leq v$, $\alpha,\gamma\in\Rmaxm$,
$\alpha\leq \gamma$, will be called {\em (idempotent) affine halfspace}.
It is called {\em archimedean}, if  $u$, $v$, $\alpha$ and $\gamma$
are positive.

\msn
For a convex set 
$C\subset \Rmaxmn$  define $V(C)\subset \R_{\max,m}^{n+1}$
to be the semimodule of vectors
of the form $(x_1\lambda,\ldots,x_n\lambda,\lambda)$ with 
$x=(x_1,\ldots,x_n)\in C$ and $\lambda\in \Rmaxm$.

\begin{theorem}
\label{maxsep3}
{\rm (Separation of convex sets)}
Let $C_1,\ldots,C_k$ be closed convex subsets of $\R_{\max,m}^n$
with empty intersection, and
assume that the intersection of the recession cones of $C_1,\ldots,C_k$
is zero.
Then, there exist affine archimedean 
halfspaces $H^{\aff}_1,\ldots,H^{\aff}_k$ which 
contain the corresponding convex sets $C_i,\, i=1,\ldots,k$ 
and have empty intersection.
\end{theorem}
\begin{proof}
From Prop.~2.16 of \cite{GK-07<LAA>},
we know that the closure of $V(C_i)$, $\overline{V(C_i)}$,
is equal to $V(C_i)\cup (\operatorname{rec}(C_i)\times \{0\})$.
Hence, the assumptions imply that the intersection of
$\overline{V(C_1)},\ldots,\overline{V(C_k)}$
is zero. By Theorem~\ref{maxsep2}, we can find archimedean
halfspaces $H_i\supset \overline{V(C_i)}$ with zero intersection.
Every $H_i$ can be written
as 
\begin{equation}
\label{e:h'i}
H_i=\{(x_1,\ldots,x_n,\mu)\mid 
u^i/x\wedge \alpha^i/\mu \geq v^i/x\wedge \gamma^i/\mu\}\cup \{0\}
\end{equation}
with $u^i\leq v^i$ and $\alpha^i\leq \gamma^i$, 
understanding that $x:=(x_1,\ldots,x_n)$.
Observe that for all $x\in C_i$, $(x,1)\in V(C_i)\subset H_i$.
We deduce that the affine archimedean halfspace
\begin{equation}
\label{hi10}
H^{\aff}_i=\{x\mid u^i/x\wedge \alpha^i \geq v^i/x\wedge \gamma^i\}
\end{equation}
contains $C_i$. 
Since the intersection 
of the halfspaces $H_i$  is zero, 
the intersection of the affine halfspaces $H^{\aff}_i$ must be empty.
\end{proof}\eproof

\msn In convex analysis, one can find an analogous separation theorem 
for several compact convex sets, see \cite{Egg:58}, pages 39-40.

\msn We now deduce an idempotent analogue of the classical
Helly's Theorem. As observed by S. Gaubert and
F. Meunier \cite{Meu-06}, there is another proof
of this theorem, which is based on the direct
idempotent analogue of Radon's argument (see \cite{Egg:58}).

\begin{theorem}
\label{helly1}
{\rm (Helly's Theorem)}
Suppose that $V_i,\ i=1,\ldots,m$
is a collection of $m\geq n$ semimodules in $\Rmaxmn$.
If each $n$ semimodules intersect nontrivially, then the whole collection
has a nontrivial intersection.
\end{theorem}
\begin{proof} 
It suffices to consider the case
where the semimodules $V_i$ are all closed. Indeed, the assumption
implies that for all $j:=(j_1,\ldots,j_n)\in \{1,\ldots,m\}^n$, we
can choose a non-zero element $z_j$ in the intersection 
$V_{j_1}\cap\cdots \cap V_{j_n}$. 
Let $V'_i$ denote the semimodule generated by the elements 
$z_j$ that belong to $V_i$. 
The collection of semimodules $V'_i$, $i=1,\ldots,m$ still has 
the property that each $n$ semimodules intersect nontrivially. Moreover,
$V'_i$ is closed, because it is finitely generated (see
e.g. Lemma~2.20 of \cite{GK-07<LAA>} or Corollary~27 of 
\cite{BSS-07}). Hence, if the conclusion
of the theorem holds for closed semimodules, 
we deduce that the whole collection 
$V'_i,\, i=1,\ldots,m$ has a nontrivial intersection, and 
since $V_i\supset V'_i$, the conclusion of 
the theorem also holds without any closure assumption.

\msn In the discussions that follow, 
the semimodules $V_i$ are all closed.
We argue by contradiction, assuming that the whole
collection has zero intersection. 
By Theorem~\ref{include}, we can also assume that the semimodules
$V_i$ are archimedean.
For some number $k<m$ every 
$k$ semimodules intersect nontrivially, but there are $k+1$
semimodules, say
$V_1,\ldots,V_{k+1}$, which have zero intersection.
By Theorem \ref{maxsep}, there exists a positive vector
$y=y^0$ and a scalar $\lambda<1$ such that
\begin{equation}
\label{projineq3}
P_{k+1}\cdots P_1 y\leq\lambda y.
\end{equation}
For all $i$ we denote
$y^i=P_i\cdots P_1 y^0$, where projectors are indexed
modulo $(k+1)$.
By the homogeneity and isotonicity of projectors, we have that 
\begin{equation}
\label{projineq4}
P_{l+k+1}\cdots P_{l+1} y^l\leq\lambda y^l
\end{equation}
for all $l=1,\ldots$. Consider the vectors
\begin{equation}
\label{vects}
z^l=P_{l+k}\cdots P_{l+1} y^l
\end{equation}
for $l=1,\ldots,k+1$. Since each $k$ semimodules intersect 
nontrivially, the vector $z^l$ must have at least one
coordinate equal to that of $y^l$, for otherwise
$y^l$ would satisfy the first condition of
Theorem \ref{maxsep}, giving a contradiction.   
As $k\geq n$, there are at least two numbers $l$ and at least one number
$i$ such that $z^l$
has the same $i$th coordinate as $y^l$. 
If we take the smallest of these
two $l$ numbers, then
\begin{equation}
\label{sleeper}
(P_{l+k+1}\cdots P_{l+1} y^l)_i=y^l_i.
\end{equation}
But this contradicts (\ref{projineq4}). Hence any
$k+1$ semimodules intersect
nontrivially, which is again a contradiction. The theorem
is proved. \end{proof}\eproof

There is also an affine version of this theorem.

\begin{theorem}
\label{helly-aff}
Suppose that $C_i,\, i=1,\ldots,m$ is a collection
of $m\geq n+1$ convex subsets of $\R_{\max,m}^n$. If each $n+1$ of these
convex sets have a nonempty intersection, then the whole collection
has a nonempty intersection.
\end{theorem}
\begin{proof}
Consider the
semimodules $V(C_1),\ldots,V(C_m)$ defined above, 
and apply Theorem~\ref{helly1} to them. 
\end{proof}\eproof

\section{Acknowledgements}

The two authors thank 
Peter Butkovi\v{c} and 
Hans Schneider
for illuminating discussions 
which have been at the origin
of the present work.
The first author also thanks Fr\'{e}d\'{e}ric Meunier
for having drawn his attention to
the max-plus analogues of Helly-type 
theorems. The second author is grateful to 
Andre\u{\i} Sobolevski\u{\i} for valuable 
ideas and discussions
concerning the analogy between convex geometry and idempotent analysis.


\end{document}